\documentclass{amsart}
\usepackage{newsymbol}
\newsymbol\forces   130D
\newsymbol\onto     1310
\newsymbol\restr    1216
\newsymbol\into     131A
\newsymbol\bale     1334  \newcommand\baledot{\mathrel{\dot\bale}}
\newsymbol\le       1336
\newsymbol\bage     133C  
\newsymbol\ge       133E
\newsymbol\lelex    1343
\newsymbol\notle    230A
\newsymbol\notge    230B
\newsymbol\emptyset 203F
\newbbbletter\baB       B
\newbbbletter\baC       C \let\complex=\baC
\newbbbletter\halfline  H
\newbbbletter\unitint   I
\newbbbletter\Miod      M
\newbbbletter\naturals  N
\newbbbletter\Poset     P
\newbbbletter\rationals Q
\newbbbletter\reals     R
\newbbbletter\baT       T
\newbbbletter\Z         Z
\newcommand  \cont      {\mathfrak{c}}
\newcommand  \Hstar     {\halfline\mskip1mu{}^*}
\newcommand  \Isubu     {\unitint_u}
\newcommand  \betaH     {\beta\halfline}
\newcommand  \gammaN    {\gamma\omega}
\newcommand  \Mstar     {\Miod^*}
\newcommand  \Nstar     {\omega^*}
\newcommand  \Rstar     {\reals^*}
\newcommand  \LL        {\mathfrak{L}}
\newcommand  \cl        {\operatorname{cl}}
\renewcommand\int       {\operatorname{int}}
\newcommand  \Fn        {\operatorname{Fn}}

\newcommand  \dom       {\operatorname{dom}}
\newcommand  \tc        {\operatorname{trcl}}
\newcommand  \inv       {^{\scriptscriptstyle\leftarrow}}
\newcommand  \Hyper[1]  {2^{#1}}
\newcommand  \seq[2][n] {\langle #2_#1\rangle_#1}
\newcommand  \card[1]   {\lvert#1\rvert}
\newcommand  \bigcard[1]{\bigl|#1\bigr|} \let\bigabs\bigcard
\newcommand  \0         {\mathbf{0}}
\newcommand  \orpr[2]   {\langle#1,#2\rangle}
\DeclareMathSymbol \calB {0}{symbols}{`B}
\DeclareMathSymbol \calC {0}{symbols}{`C}
\DeclareMathSymbol \calD {0}{symbols}{`D}
\DeclareMathSymbol \calF {0}{symbols}{`F}
\DeclareMathSymbol \calI {0}{symbols}{`I}
\DeclareMathSymbol \calL {0}{symbols}{`L}
\DeclareMathSymbol \calM {0}{symbols}{`M}
\DeclareMathSymbol \pow  {0}{symbols}{`P}
\DeclareMathSymbol \calR {0}{symbols}{`R}
\let\axiom\mathsf
\newcommand\CH {\axiom{CH}}
\newcommand\GCH{\axiom{GCH}}
\newcommand\MA {\axiom{MA}}
\newcommand\PFA{\axiom{PFA}}
\newcommand\ZFC{\axiom{ZFC}}
\usepackage{boldmath}
\newboldletter\bH H1
\newboldletter\bN N1
\DeclareMathSymbol\land {3}{symbols}{"5E}
\DeclareMathSymbol\lor  {3}{symbols}{"5F}
\DeclareMathSymbol\diag {1}{symbols}{"34}
\let\phi\varphi          \newcommand\phidot{\dot\phi}
\let\liminv\varprojlim

\newtheorem{introlemma}{Lemma}[section]

\newtheorem{theorem}{Theorem}[section]
\renewcommand{\thetheorem}{\arabic{section}.\arabic{theorem}}
\newtheorem{lemma}[theorem]{Lemma}
\newtheorem{falselemma}[theorem]{False lemma}
\theoremstyle{definition}
\newtheorem{example}[theorem]{Example}
\theoremstyle{remark}
\newtheorem{remark}[theorem]{Remark}
\newtheorem{question}[theorem]{Question}
\newtheorem{problem}[theorem]{Problem}
\def\iff/{if\kern0pt f}

\makeatletter
\def\cases{\left\{\,\vcenter\bgroup\normalbaselines\m@th\def\\{\cr}
           \ialign\bgroup$##\hfil$&\quad##\hfil\crcr}
\def\endcases{\crcr\egroup\egroup\right.}
\def\eqalign{\null\,
             \vcenter\bgroup\openup\jot\m@th\def\\{\cr}
                     \ialign\bgroup\strut\hfil$\displaystyle{##}$&
                                        $\displaystyle{{}##}$\hfil\crcr}
\def\endeqalign{\crcr\egroup\egroup\,}
\makeatother
\let\cyr\relax
\newcommand\cprime{{\mathsurround0pt$'$}}
\begin{document}

\title{A Universal Continuum of Weight~$\aleph$}

\author{Alan Dow}
\address{Department of Mathematics and Statistics\\
         York University\\
         4700 Keele Street\\
         Toronto, Ontario\\
         Canada M3J 1P3}
\email{dowa@yorku.ca}

\author[Klaas Pieter Hart]{Klaas Pieter Hart*}
\thanks{*The research of the second author was supported by the Netherlands
        Organization for Scientific Research (NWO) --- Grant~R\,61-322.}
\address{Faculty of Technical Mathematics and Informatics\\
         TU Delft\\
         Postbus 5031\\
         2600~GA~~Delft\\
         the Netherlands}
\email{k.p.hart@twi.tudelft.nl}

\begin{abstract}
We prove that every continuum of weight~$\aleph_1$ is a continuous image
of the \v{C}ech-Stone-remainder~$\Rstar$ of the real~line.
It follows that under~$\CH$ the remainder of the half~line
$[0,\infty)$ is universal among the continua of weight~$\cont$
--- universal in the `mapping onto' sense.

We complement this result by showing that
\ 1)~under~$\MA$ every continuum of weight less than~$\cont$ is a continuous
image of\/~$\Rstar$,
\ 2)~in the Cohen model the long segment of length~$\omega_2+1$ is not a
continuous image of\/~$\Rstar$, and
\ 3)~$\PFA$ implies that $\Isubu$ is not a continuous image of\/~$\Rstar$,
     whenever $u$~is a $\cont$-saturated ultrafilter.

We also show that a universal continuum can be gotten from a $\cont$-saturated
ultrafilter on~$\omega$ and that it is consistent that there is no universal
continuum of weight~$\cont$.
\end{abstract}



\maketitle

\section*{Introduction}

\renewcommand\thetheorem{\arabic{theorem}}

In \cite{Parovicenko63} Parovi\v{c}enko proved that every compact
Hausdorff space of weight~$\aleph_1$ is a continuous image of the
space~$\Nstar$
(the \v{C}ech-Stone remainder of the discrete space~$\omega$).
It follows that under~$\CH$ the space~$\Nstar$ is a universal compact
space of weight~$\cont$, universal in the sense of onto mappings
rather than embeddings.

The purpose of this paper is to prove a similar result for $\Hstar$,
the \v{C}ech-Stone remainder of the half~line $\halfline=[0,\infty)$.
As $\Hstar$ is compact and connected (a \emph{continuum}) the following
theorem is the proper analogue of Parovi\v{c}enko's result.

\begin{theorem}\label{main.thm}
Every continuum of weight\/~$\aleph_1$ is a continuous image of\/~$\Hstar$.
\end{theorem}

We present proofs of our theorem using each of three quite different
approaches to Parovi\v{c}enko's theorem from the literature that have
proven valuable.

It is clear that new ideas are needed to prove Theorem~\ref{main.thm}
since the fact that $\Nstar$~is zero-dimensional appears to be so
essential in each of the proofs of Parovi\v{c}enko's theorem.

We now give brief outlines of the proofs of Parovi\v{c}enko's theorem
and we indicate how our proofs parallel these arguments; more details,
including definitions, will appear later in the paper.
We shall finish this introduction with some remarks on universal
spaces.

\subsection*{Two algebraic proofs}

The first two proofs of Parovi\v{c}enko's theorem use Stone's duality,
see~\cite{StoneMH37a}, between compact zero-di\-men\-sional spaces and
Boolean algebras.
One first uses a theorem of Alexandroff, from~\cite{Alexandroff36},
to find, given the compact space~$X$ of weight~$\aleph_1$, a compact
zero-di\-men\-sional space~$Y$ of weight~$\aleph_1$ and a continuous
map from~$Y$ onto~$X$, thus reducing the problem to the case of a
zero-di\-men\-sional range space.
Stone Duality then tells us that finding a continuous map from~$\Nstar$
onto~$Y$ is equivalent to finding an embedding of the Boolean algebra~$\baC$
of clopen subsets of\/~$Y$ into the Boolean algebra~$\baB$ of clopen subsets
of\/~$\Nstar$.

This gives rise to two proofs of Parovi\v{c}enko's theorem.
In the first, one constructs the desired embedding by transfinite
induction; this is Parovi\v{c}enko's original proof.
The second proof recognizes that Parovi\v{c}enko's theorem is in fact
a special case of a well-known result about $\aleph_1$-saturated models.

Our parallels of these proofs are naturally based on Wallman's generalization
of Stone's duality to the class of distributive lattices and amount to
embedding, for the given continuum, a suitable sublattice of its lattice of
closed sets into the lattice of closed subsets of~$\Hstar$.
The key new idea is to extract the hidden model-theoretic aspects of the
aforementioned proofs.

\subsection*{A topological proof}

The third proof is of a more topological nature, it appears
in~\cite{BlaszczykSzymanski80}.
It starts by writing the compact space~$X$ as an inverse
limit, $\liminv X_\alpha$, of an $\omega_1$-sequence of compact metrizable
spaces.
The continuous map from~$\Nstar$ onto~$X$ is obtained as the limit of a
family of continuous onto maps~$f_\alpha:\Nstar\onto X_\alpha$.
These maps are constructed by transfinite induction and using the following
lemma.

\begin{introlemma}\label{lemma.intro}
Let $M$ and $N$ be compact metrizable spaces and let $f:\Nstar\onto N$
and $g:M\onto N$ be continuous surjections.
Then there is a continuous surjection $h:\Nstar\onto M$ such that
$g\circ h=f$.
\end{introlemma}

Our third proof will be along these lines and the main problem is that
the above lemma is false if we simply replace $\Nstar$ by $\Hstar$ and
assume that $M$ and $N$ are continua --- see Example~\ref{example.wish}.
Intriguingly, if $g$~is induced by an appropriate elementary embedding
then the analogue of Lemma~\ref{lemma.intro} will hold; indeed,
this is the key new feature of our proofs: the introduction of elementarity.

\subsection*{Universal spaces}

As with Parovi\v{c}enko's theorem, Theorem~\ref{main.thm} becomes
a result about universal spaces if one assumes
the Continuum Hypothesis ($\CH$).

\begin{theorem}[$\CH$]
The space\/ $\Hstar$~is a universal continuum of weight~$\cont$.
\end{theorem}

It is of interest to note that in the class of compact metrizable spaces
the situation is somewhat different from what we just found for
spaces of weight~$\cont$.
On the one hand the Cantor set \emph{is} a universal compact metrizable space;
this is a theorem of Alexandroff~\cite{Alexandroff27} and
Hausdorff~\cite{Hausdorff27}.
On the other hand, in~\cite{Waraszkiewicz34} Waraszkiewicz constructed
a family of plane continua such that no metric continuum maps onto all
of them, hence there is no universal metric continuum.

Part of the balance is restored by the following theorem,
due to Aarts and Van Emde Boas \cite{AartsvanEmdeBoas67};
as we shall need this theorem and its proof later we provide a short
argument.

\begin{theorem}\label{thm.jan}
The space\/ $\Hstar$ maps onto every metric continuum.
\end{theorem}

\begin{proof}
Consider a metric continuum~$M$ and assume it is embedded into
the Hilbert cube~$Q=[0,1]^\infty$.
Choose a countable dense subset~$A$ of\/~$M$ and enumerate it
as~$\{a_n:n\in\omega\}$.
Next choose, for every~$n$, a finite sequence of points
~$a_n=a_{n,0}$, $a_{n,1}$, \dots,~$a_{n,k_n}=a_{n+1}$
such that $d(a_{n,i},a_{n,i+1})<2^{-n}$ for all~$i$ --- this uses the
connectivity of\/~$M$.
Finally, let~$e$ be the map from~$\halfline$
to~$(0,1]\times Q$ with first coordinate~$e_1(t)=2^{-t}$
and whose second coordinate satisfies $e_2(n+\frac i{k_n})=a_{n,i}$
for all~$n$ and~$i$ and is (piecewise) linear otherwise.

It is clear that $e$ is an embedding and one readily checks that
$\cl e[\halfline]=e[\halfline]\cup\bigl(\{0\}\times M\bigr)$;
the \v{C}ech-Stone extension~$\beta e$ of\/~$e$ maps~$\Hstar$ onto~$M$.
\end{proof}

\subsection*{Limitations and extensions}

As with Parovi\v{c}enko's results ours are limited to~$\aleph_1$:
we show by examples similar to the ones used for~$\Nstar$ that
not necessarily every continuum of weight~$\aleph_2$ is a continuous image
of~$\Hstar$ and that Martin's Axiom is not strong enough to guarantee
that every continuum of weight~$\cont$ is a continuous image of~$\Hstar$.

On the other hand, in Section~\ref{sec.more.univ.continua} we exhibit
a few more continua with the same behaviour as~$\Hstar$ and toward the end
of Section~\ref{sec.ext.and.lim} we show that from a $\cont$-saturated
ultrafilter one can construct a universal continuum of weight~$\cont$.

The ultimate limitation occurs at the end of the paper, where we show it
consistent that no universal continuum of weight~$\cont$ exists.

\smallskip
We thank the referee for a detailed and thoughtful report; it enabled us to
improve both content and presentation of this paper.

\renewcommand\thetheorem{\arabic{section}.\arabic{theorem}}

\section{Algebraic Preliminaries}
\label{sec.algebra}

We already indicated that two of our proofs will involve the lattice
of closed sets of our continua and that the basis for this is
Wallman's generalization, to the class of distributive lattices,
of Stone's representation theorem for Boolean algebras.
Wallman's representation theorem is as follows.

\begin{theorem}[Wallman \cite{Wallman38}]
If $L$ is a distributive lattice then there is a compact $T_1$-space~$X$
with a base for its closed sets that is a homomorphic image of\/~$L$.
The homomorphism is an isomorphism \iff/ $L$ is \emph{disjunctive},
which means: if $a\notle b$ then there is $c\in L$ such that $c\le a$
and $c\wedge b=0$.
\end{theorem}

The representing space~$X$ is Hausdorff \iff/ the lattice~$L$ is
\emph{normal}; this means that given $a,b\in L$ with $a\wedge b=0$
there are $c,d\in L$ such that $c\vee d=1$, \ $a\wedge d=0$ and $b\wedge c=0$.
Note that this mimics the formulation of normality of topological spaces
in terms of closed sets only.

Also, if one starts out with a compact Hausdorff space~$X$ and a
base~$\calC$ for its closed sets that is closed under finite unions
and intersections then this base~$\calC$ is a normal lattice and its
representing space is~$X$.

\smallskip

The following theorem shows how one can create onto mappings
from maps between lattices.
In it we use $\Hyper{X}$ to denote the family of closed subsets of the
space~$X$.

\begin{theorem}\label{thm.how.to.map.onto}
Let $X$ and $Y$ be compact Hausdorff spaces and let $\calC$~be a base for the
closed subsets of\/~$Y$ that is closed under finite unions and finite
intersections.
Then $Y$~is a continuous image of\/~$X$ if and only if there is a map
$\phi:\calC\to\Hyper{X}$ such that
\begin{enumerate}
\item $\phi(\emptyset)=\emptyset$ and\label{cond.i}
      if $F\neq\emptyset$ then $\phi(F)\neq\emptyset$;
\item if $F\cup G=Y$ then $\phi(F)\cup\phi(G)=X$; and\label{cond.ii}
\item if $F_1\cap\cdots\cap F_n=\emptyset$\label{cond.iii}
      then $\phi(F_1)\cap\cdots\cap\phi(F_n)=\emptyset$.
\end{enumerate}
\end{theorem}

\begin{proof}
Necessity is easy: given a continuous onto map $f:X\to Y$ let
$\phi(F)=f\inv[F]$.
Note that $\phi$ is in fact a lattice-embedding.

To prove sufficiency let $\phi:\calC\to\Hyper{X}$ be given and consider
for each $x\in X$ the family
$\calF_x=\bigl\{F\in\calC:x\in\phi(F)\bigr\}$.
We claim that $\bigcap\calF_x$ consists of exactly one point.
Indeed, by Condition~\ref{cond.iii} the family~$\calF_x$ has the finite
intersection property, so that $\bigcap\calF_x$ is nonempty.
Next assume that $y_1\neq y_2$ in~$Y$ and take $F,G\in\calC$ such that
$F\cup G=Y$,\ $y_1\notin F$ and $y_2\notin G$.
Then, by Condition~\ref{cond.ii}, either $x\in\phi(F)$ and so
$y_1\notin\bigcap\calF_x$ or $x\in\phi(G)$ and so $y_2\notin\bigcap\calF_x$.

We define $f(x)$ to be the unique point in~$\bigcap\calF_x$.

To demonstrate that $f$ is continuous and onto we show that for every closed
subset~$F$ of\/~$Y$ we have
$$
f\inv[F]=\bigcap\bigl\{\phi(G):G\in\calC, F\subseteq\int G\bigr\}.\eqno(*)
$$
This will show that preimages of closed sets are closed and that every fiber
$f\inv(y)$~is nonempty.

We first check that the family on the right-hand side has the finite
intersection property.
Even though $F$ and the complement~$K$ of~$\bigcap_i\int G_i$ need not
belong to~$\calC$ we can still find~$G$ and~$H$ in~$\calC$ such that
$G\cap K=H\cap F=\emptyset$ and $H\cup G=Y$.
Indeed, apply compactness and the fact that $\calC$~is a lattice to
find~$C$ in~$\calC$ such that $F\subseteq C\subseteq\bigcap_i\int G_i$ and
then $D\in\calC$ with $K\subseteq D$ and $D\cap C=\emptyset$;
then apply normality of~$\calC$ to~$C$ and~$D$.
Once we have $G$ and~$H$ we see that for each~$i$ we also have~$H\cup G_i=Y$
and so $\phi(H)\cup\phi(G_i)=X$; combined with
$\phi(G)\cap\phi(H)=\emptyset$ this gives
$\phi(G)\subseteq\bigcap_i\phi(G_i)$.

To verify~$(*)$, first let $x\in X\setminus f\inv[F]$.
As above we find $G$ and~$H$ in~$\calC$ such that
$f(x)\notin G$, $H\cup G=Y$ and $H\cap F=\emptyset$.
The first property gives us~$x\notin\phi(G)$; the other two imply
that $F\subseteq\int G$.

Second, if $F\subseteq\int G$ then one can find $H\in\calC$ such that
$H\cup G=X$ and $H\cap F=\emptyset$.
It follows that if $x\notin\phi(G)$ we have $x\in\phi(H)$, hence $f(x)\in H$
and so $f(x)\notin F$.
\end{proof}

An obvious corollary of this theorem is that $Y$ is a continuous image
of\/~$X$ \iff/ there is an embedding of the lattice~$\calC$ into the
lattice~$\Hyper{X}$.

Theorem~\ref{thm.how.to.map.onto} is formulated so as to be applied in
the following way:
given the continuum~$K$ we shall take a base~$\calC$ for its closed sets
and define a map as in the theorem into the canonical base~$\calL$ for
the closed sets of\/~$\Hstar$, rather than into the
family~$\Hyper{\Hstar}$.
We do not know whether it is actually possible to get a lattice-embedding
of\/~$\calC$ into~$\calL$.

The canonical base in question is just
$$
\calL=\{A^*:\text{$A$~is closed in~$\halfline$}\,\}.
$$
Here, as is common, $A^*$ abbreviates $\cl A\cap\Hstar$.
Notice that this is exactly analogous to the collection of clopen
subsets of $\omega^*$ but that it also points up another quite stark
difference since unlike in the zero-dimensional case there doesn't
appear to be a suitable internally defined base to choose.
The base consisting of all zero sets is worth considering but it lacks the
model theoretic saturation properties that we will require.

The following lemma will be the first step in the constructions of both
embeddings.

\begin{lemma}\label{lemma.first.step}
Let $K$ be a metric continuum and let  $x\in K$.
There is a map~$\phi$ from~$\Hyper{K}$ to~$\Hyper{\halfline}$ such that
\begin{enumerate}
\item $\phi(\emptyset)=\emptyset$ and $\phi(K)=\halfline$;\label{fs.i}
\item $\phi(F\cup G)=\phi(F)\cup\phi(G)$;\label{fs.ii}
\item if $F_1\cap\cdots\cap F_n=\emptyset$\label{fs.iii}
      then $\phi(F_1)\cap\cdots\cap\phi(F_n)$ is compact; and
\item $\naturals\subseteq\phi\bigl(\{x\}\bigr)$.\label{fs.iv}
\end{enumerate}
In addition, if some countable family~$\calC$ of nonempty closed subsets
of\/~$K$ is given in advance then we can arrange that for every~$F$ in~$\calC$
the set~$\phi(F)$ is not compact.
\end{lemma}

\begin{proof}
As proved in Theorem~\ref{thm.jan} there is a map from~$\Hstar$ onto~$K$.

The proof given in~\cite{AartsvanEmdeBoas67} is flexible enough to allow us
to ensure that the embedding~$e$ of\/~$\halfline$ into~$(0,1]\times Q$ is such
that $e(n)=\langle2^{-n},x\rangle$ for every~$n\in\naturals$
and that for every element~$y$ of some countable set~$C$ the set
$\{t:e_1(t)=y\}$ is cofinal in~$\halfline$ --- it is also easy to change the
description of~$e$ in the proof we gave to produce another~$e$ with the
desired properties.
In our case we let $C$ be a countable subset of\/~$K$ that meets every
element of the family~$\calC$.

We now identify $K$ and $\{0\}\times K$ and define a map
$\psi:\Hyper{K}\to\Hyper{\unitint\times Q}$ by
$$
\psi(F)=\bigl\{y\in:\unitint\times Q:d(y,F)\le d(y,K\setminus F)\bigr\}.
$$
In \cite[\S\,21\,XI]{Kuratowski66} it is shown that for all $F$ and $G$
we have
\begin{itemize}
\item $\psi(F\cup G)=\psi(F)\cup\psi(G)$;
\item $\psi(K)=\unitint\times Q$ and $\psi(\emptyset)=\emptyset$ ---
      by the fact that $d(y,\emptyset)=\infty$ for all~$y$; and
\item $\psi(F)\cap K=F$.
\end{itemize}
Note that for every~$y\in K$ we have
$d\bigl(\langle t,y\rangle,\{y\}\bigr)
 =d\bigl(\langle t,y\rangle,K\setminus\{y\}\bigr)=t$ and hence
$\unitint\times\{y\}\subseteq\psi\bigl(\{y\}\bigr)$.

Now define $\phi(F)=e\inv\bigl[\psi(F)\bigr]$ or rather, after identifying
$\halfline$ and~$e[\halfline]$, set $\phi(F)=\psi(F)\cap e[\halfline]$.
All desired properties are easily verified: \ref{fs.i} and \ref{fs.ii}
are immediate;
to see that \ref{fs.iii}~holds note that
if $F_1\cap\cdots\cap F_n=\emptyset$
then $\cl\phi(F_1)\cap\cdots\cap\cl\phi(F_n)\cap K=\emptyset$,
so that $\cl\phi(F_1)\cap\cdots\cap\cl\phi(F_n)$ is a compact subset
of\/~$\halfline$.
That \ref{fs.iv}~holds follows from the way we chose the values~$e(n)$ for
$n\in\naturals$.

Finally, if $F\in\calC$ and $y\in C\cap F$ then
the cofinal set $\{t:\pi\bigl(e(t)\bigr)=y\}$ is a subset of\/~$\phi(F)$
so that $\phi(F)$~is not compact.
\end{proof}

\begin{remark}\label{rem.lattice.hom}
If $\phi$ is as in Lemma~\ref{lemma.first.step} then the
map~$\psi:\calC\to\calL$ defined by
$$
\psi:C\mapsto \phi(C)\mapsto \phi(C)^*
$$
is as in Theorem~\ref{thm.how.to.map.onto}.
Observe that $\psi$ is even a $\cup$-homomorphism; to get a lattice
homomorphism the map~$\phi$ would have to satisfy~$3'$ instead of\/~$3$,
where $3'$~reads:
\begin{enumerate}
\item[$3'$.] the symmetric difference of $\phi(F\cap G)$ and
             $\phi(F)\cap\phi(G)$ is bounded in~$\halfline$.
\end{enumerate}
We do not know how to achieve~$3'$ and at this point there is no extra
benefit to be had from this condition so we leave it be.
\end{remark}

\section{Two algebraic proofs}

We shall now show how to construct, given a continuum $K$, a map~$\phi$
from a base for the closed sets of~$K$ into the base~$\calL$
as in Theorem~\ref{thm.how.to.map.onto}.

\subsection*{A proof using Model Theory}

Our plan is to find the map~$\phi$ promised above by an appeal to
some machinery from Model Theory.
The first step would be to show that $\calL$~is an $\aleph_1$-saturated
lattice and hence an $\aleph_2$-universal structure; this latter notion
says that every structure of size~$\aleph_1$ that is
\emph{elementarily equivalent} to~$\calL$ is embeddable into~$\calL$.
The second step would then be to show that every lattice of
size~$\aleph_1$ is embeddable into a lattice of size~$\aleph_1$ that
itself is elementarily equivalent to~$\calL$.

There are two problems with this approach:
1)~we were not able to show that $\calL$~is
$\aleph_1$-saturated, and
2)~Lemma~\ref{lemma.first.step} does not give a lattice-embedding.
We shall deal with these problems in turn.
But first we give some definitions from Model Theory.
We shall try to illustrate these definitions by means of linear orderings.

\subsubsection*{Saturation and universality}

Our basic reference for Model Theory is Hodges' book \cite{Hodges93},
in particular Chapter~10.
A structure (e.g., a field, a group, an ordered set, a lattice)
is said to be $\aleph_1$-saturated if, loosely speaking,
every countable consistent set of equations has a solution, where
a set of equations is consistent if every finite subsystem has a solution.
Thus, the ordered set of the reals is \emph{not} $\aleph_1$-saturated because
the following countable system of equations, though consistent, does not have a
solution: $0<x$ together with $x<\frac1n$ ($n\in\naturals$).
On the other hand, any ultrapower~$\reals^\omega_u$ of~$\reals$ is
$\aleph_1$-saturated as an ordered set --- see \cite[Theorem~9.5.4]{Hodges93}.
Such an ultrapower is obtained by taking the power~$\reals^\omega$,
an ultrafilter~$u$ on~$\omega$ and identifying points~$x$ and~$y$
if $\{n:x_n=y_n\}$ belongs to~$u$.
The ordering~$<$ is defined in the obvious way:
$x<y$ \iff/ $\{n:x_n<y_n\}$ belongs to~$u$.
It is relatively easy to show that this gives an $\aleph_1$-saturated
ordering; given a countable consistent set of equations $x<a_i$ and $x>b_i$
($i\in\omega$) one has to produce a single~$x$ that satisfies them all;
the desired~$x$ can be constructed by a straightforward diagonalization.

Two structures are \emph{elementarily equivalent} if they satisfy the same
sentences (formulas without free variables); for example
any two dense linearly ordered sets without end~points are elementarily
equivalent as linear orders, see~\cite[Section~2.7]{Hodges93}.
From this, and the fact that $\aleph_1$-saturated structures are
$\aleph_2$-universal, it follows at once that any ultrapower~$\reals^\omega_u$
contains an isomorphic copy of every $\aleph_1$-sized dense linear order without
end~points --- a result that can also be established directly by a
straightforward transfinite recursion.
It also follows that $\reals^\omega_u$ contains an isomorphic copy of
\emph{every} $\aleph_1$-sized linear order: simply make it dense by inserting a
copy of the rationals between any pair of neighbours and attach copies of the
rationals at the beginning and the end to get rid of possible end~points; the
resulting ordered set is still of cardinality~$\aleph_1$ and can therefore be
embedded into~$\reals^\omega_u$.

This then is how our proof shall run: we identify a certain
$\aleph_1$-saturated structure that serves as a base for the closed sets
of some quotient of~$\Hstar$ and we show how to expand every $\aleph_1$-sized
structure into one of the same size but elementarily equivalent to
our $\aleph_1$-saturated one.

\subsubsection*{An $\aleph_1$-saturated structure}

As mentioned above we do not know whether $\calL$~is
$\aleph_1$-saturated.
We can however find an $\aleph_1$-saturated sublattice:
$$
\calL'=\{A^*:\text{$A$~is closed in~$\halfline$, and
           $\naturals\subseteq A$ or $\naturals\cap A=\emptyset$}\,\}.
$$
This lattice is a base for the closed sets of the space~$\bH$, obtained
from~$\Hstar$ by identifying~$\Nstar$ to a point.

To see that $\calL'$ is $\aleph_1$-saturated we introduce another space,
namely $\Miod=\omega\times\unitint$, where $\unitint$~denotes the unit
interval.
The canonical base~$\calM$ for the closed sets of~$\Mstar$ is naturally
isomorphic to the reduced power $(2^\unitint)^\omega$ modulo the
cofinite filter.
It is well-known that this structure is $\aleph_1$-saturated
--- see~\cite{JonssonOlin68}.
The following substructure~$\calM'$ is $\aleph_1$-saturated as well:
$$
\calM'=\{A^*:\text{$A$~is closed in~$\Miod$, and
             $\bN\subseteq A$ or
             $\bN\cap A=\emptyset$}\,\},
$$
where $\bN=\{0,1\}\times\omega$.
Indeed, consider a countable set $T$ of equations with constants
from~$\calM'$ such that every finite subset has a solution in~$\calM'$.
We can then add either either $\bN\subseteq x$ or $\bN\cap x=\emptyset$
to~$T$ without losing consistency.
Any element of~$\calM$ that satisfies the expanded~$T$ will automatically
belong to~$\calM'$.

We claim that $\calL'$ and $\calM'$ are isomorphic.
To see this consider the map $q:\Miod\to\halfline$ defined by
$q(n,x)=n+x$.
The \v{C}ech-Stone extension of~$q$ maps $\Miod^*$ onto~$\Hstar$ and it
is readily verified that $L\mapsto q^{-1}[L]$ is an isomorphism
between~$\calL'$ and~$\calM'$.
(In topological language: the space~$\bH$ is also obtained from~$\Mstar$
by identifying~$\bN^*$ to a point.)

\subsubsection*{A new language}

The last point that we have to address is that Lemma~\ref{lemma.first.step}
does not provide a lattice embedding but rather a map that only partially
preserves unions and intersections.
This is where Theorem~\ref{thm.how.to.map.onto} comes in:
we do not need a full lattice embedding but only a map that preserves certain
identities.
We abbreviate these identities as follows:
$$
\begin{aligned}
J(x,y) &\quad\equiv\quad x\vee y=1 \\
M_n(x_1,\ldots,x_n) &\quad\equiv\quad x_1\wedge\cdots\wedge x_n=0.
\end{aligned}
$$
We can restate the conclusion in Theorem~\ref{thm.how.to.map.onto} in the
following manner: $Y$~is a continuous image of~$X$ if and only if there is
an~$\LL$-homomorphism from~$\calC$ to~$\Hyper{X}$, where $\LL$~is the language
that has $J$ and the~$M_n$ as its predicates \emph{and} where $J$ and the~$M_n$
are interpreted as above.

Note that by considering a lattice with $0$ and~$1$ as an $\LL$-structure
we do not have to mention $0$ and~$1$ anymore; they are implicit in the
predicates.
For example we could define a normal $\LL$-structure to be one in which
$M_2(a,b)$ implies
$(\exists c,d)\bigl(M_2(a,d)\land M_2(c,b)\land J(c,d)\bigr)$.
Then a lattice is normal \iff/ it is normal as an $\LL$-structure.

\subsubsection*{The proof}

Let $\calC$ be a base of size~$\aleph_1$ for the closed sets of the
continuum~$K$.
We want to find an $\LL$-structure~$\calD$ of size~$\aleph_1$
that contains~$\calC$ and that is elementarily equivalent to~$\calL'$.
To this end we consider the diagram of\/~$\calC$; that is, we add the elements
of\/~$\calC$ to our language~$\LL$ and we consider the set~$D_\calC$ of all
atomic sentences from this expanded language that hold in~$\calC$.
For example if $a\cap b=\emptyset$ and $c\cup d=K$ then
$M_2(a,b)\land J(c,d)$ belongs to~$D_\calC$.

To $D_\calC$ we add the theory~$T_{\calL'}$ of\/~$\calL'$,
to get a theory~$T_\calC$.
Let $\calC'$ be any countable subset of~$\calC$ and assume, without loss of
generality, that $\calC'$~is a normal sublattice of~$\calC$.
The Wallman space~$X$ of~$\calC'$ is metrizable, because $\calC'$~is countable,
and connected because it is a continuous image of~$K$.
We may now apply Lemma~\ref{lemma.first.step} to obtain an $\LL$-embedding
of~$\calC'$ into~$\calL'$; indeed, Condition~4 says that $\Nstar$~will be mapped
onto a fixed point~$x$ of~$X$.

This shows that, for every countable subset~$\calC'$ of~$\calC$ the union
of $D_{\calC'}$ and~$T_{\calL'}$ is consistent and so, by the compactness
theorem, the theory~$T_\calC$ is consistent.
Let $\calD$ be a model for~$T_\calC$ of size~$\aleph_1$.
This model is as required: it satisfies the same sentences
as~$\calL'$ and it contains a copy of\/~$\calC$, to wit the set of
interpretations of the constants from~$\calC$.

\subsubsection*{An extension}
It is clear from the above reasoning and Lemma~\ref{lemma.alan} below
that the following theorem holds.

\begin{theorem}\label{thm.alan}
If $X$ is a continuum with an $\aleph_1$-saturated base for its closed
sets that maps onto every metric continuum
then it maps onto every continuum of weight\/~$\aleph_1$.
\end{theorem}

\begin{lemma}\label{lemma.alan}
Suppose that $X$ is a compact space which has a base $\calC$ for the
closed sets that is $\aleph_1$-saturated as a lattice.
Then, if $X$ maps onto a metric space~$M$, there is an embedding of a
base $\calD$ for the closed subsets of $M$ into $\calC$ such as
in Theorem~\ref{thm.how.to.map.onto}, i.e., an $\LL$-embedding.
\end{lemma}

\begin{proof}
Simply fix a map $f$ from $X$ onto $M$ and let $\calD=\{D_n : n\in\omega\}$
be a base for the closed subsets of $M$ such that each $D_n$ is a regular
closed subset of\/~$M$
(that such a base, which in addition is closed under finite unions and
intersections, exists is established in~\cite{vanDouwen81d}).
We inductively choose $C_n\in\calC$ so that $f[C_n]=D_n$.
To begin we may assume that $D_0=\emptyset $ and $D_1=M$, so, of course,
we let $C_0=\emptyset$ and $C_1=X$
(note that since $\calC$ is a lattice we may assume that $\emptyset\in\calC$).
Suppose that we have chosen~$C_{n-1}$ and we must choose~$C_n$.
Since we plan to guarantee that $f[C_m]=D_m$ for each~$m$ we must have
$\bigcap\{C_k:k\in F\}=\emptyset$ for each $F\in [\omega]^{<\omega}$
such that $\bigcap\{D_k:k\in F\}=\emptyset$ and that $C_m$~will be
empty whenever $D_m$~is empty.
Therefore our only concern is that $C_k\cup C_n$ should be~$X$ whenever
$D_k\cup D_n$ is equal to~$M$.
To accomplish this we first prove that, for each closed $K\subseteq M$, there is
a $C\in\calC$ such that $f\inv\bigl[\int(K)\bigr]\subseteq C$ and $f[C]=K$.
Indeed, let $\{F_i:i\in \omega\}$ be a family of closed subsets of\/~$M$
whose union is the interior of $K$.
Also let $\{U_i: i\in\omega\}$ be a neighbourhood base for~$K$ in~$M$.
For each $i$, there are $J_i$ and $Y_i$ in~ $\calC$ such that
$f\inv[F_i]\subseteq J_i$,
\ $J_i\subseteq f\inv\bigl[\int(K)\bigr]$,
\ $f\inv[K]\subseteq Y_i$,
and $Y_i\subseteq f\inv[U_i]$.
The countable set of equations $\{J_i\subseteq x, x\subseteq Y_i:i\in\omega\}$
is clearly consistent so, because $\calC$~is $\aleph_1$-saturated,
we can choose $L\in\calC$ that contains each~$J_i$ and is contained in
each~$Y_i$.
Clearly $L$ has the desired properties.

Now for each $k<n$ such that $D_k\cup D_n=M$, let $J_k\in\calC$ be such
that $f\inv\bigl[\int(M\setminus D_k)\bigr]\subseteq J_k$ and
$f[J_k]=M\setminus\int D_k$ (recall that $D_k$ is regular closed).
If $D_k\cup D_n$ (for $k<n$) is not equal to~$M$, then let $J_k$ be the
empty set.
Also let $J_n\in\calC$ be any set such that $f(J_n)=D_n$.
Now define $C_n$ to be the union of $J_n$ with all $J_k$ for $k<n$.
Clearly $C_n\cup C_k=X$ for all $k<n$ such that $D_n\cup D_k=M$, so it
suffices to show that $f[C_n]=D_n$.
Since $J_n\subseteq C_n$, it follows that $D_n\subseteq f[C_n]$.
Let $k<n$ and assume that $D_k\cup D_n=M$.
Then $M\setminus\int D_k \subseteq D_n$
(again using that $D_k$~is regular closed).
Since $f[J_k]=M\setminus \int D_k$, it follows that $f[C_n]\subseteq D_n$.
\end{proof}

\subsection*{A direct construction}

Rather than relying on general model-theoretic results we can also give
a direct construction of an $\LL$-embedding of a base for the closed sets
of a continuum~$K$ into the lattice~$\calL$.
We do this in order to point out where the zero-dimensionality in
Parovi\v{c}enko's theorem hides an important model theoretic step.

Let $K$ be a continuum of weight~$\aleph_1$ and fix a base~$\calC$ for the
closed sets of\/~$K$ of size~$\aleph_1$ and closed under finite unions and
intersections.
We also assume $\{x\}\in\calC$ for some fixed $x\in K$.

To construct an $\LL$-embedding~$\phi$ of\/~$\calC$ into the base~$\calL$ for
the closed sets of\/~$\Hstar$ we construct a
map $\psi:\calC\to\Hyper{\halfline}$ that satisfies
\begin{enumerate}
\item if $F\neq\emptyset$ then $\psi(F)$ is not compact;
\item if $F\cup G=K$
      then there is a $t$ in~$\halfline$ such that
      $[t,\infty)\subseteq\psi(F)\cup\psi(G)$; and
\item if $F_1\cap\cdots\cap F_n=\emptyset$
      then $\psi(F_1)\cap\cdots\cap\psi(F_n)$ is compact.
\end{enumerate}
As in Remark~\ref{rem.lattice.hom} we then define
$\phi(F)=\psi(F)^*$.

We found it necessary to use elementary substructures to guide the induction;
as a consequence we shall build the embedding countably many steps at a time,
rather than one.

\subsubsection*{Elementarity}
Consider two structures $S$ and~$L$ for the same language (two fields,
two groups, two ordered sets) and assume that $S$ is a substructure of~$L$.
We call $S$ an \emph{elementary} substructure of~$L$ if, loosely speaking,
every equation with constants in~$S$ that has a solution in~$L$ already has
a solution in~$S$.
Thus, the field~$\rationals$ of rational numbers is \emph{not} an elementary
substructure of~$\reals$ --- consider the equation $x^2=2$ ---
but the algebraic numbers do form an elementary subfield of~$\complex$ ---
this may be gleaned from~\cite[Theorem~A.5.1]{Hodges93}.
As one can see from our applications below, an `equation' can be a system of
equations and a `solution' can be an $n$-tuple.
By way of example we show that an elementary sublattice~$S$ of a normal
lattice~$L$ is normal as well:
take disjoint elements $a$ and~$b$ of~$S$ and consider the system
$\bigl\{M_2(a,y), M_2(b,x), J(x,y)\bigr\}$;
because $L$~is normal this system has a solution in~$L$.
Elementarity guarantees that it must have a solution in~$S$ as well.

Elementary substructures abound because of the L\"owenheim-Skolem
theorem, which says that every subset~$A$ of a
structure~$L$ can be enlarged to an elementary substructure~$S$;
moreover if the language in question is countable then one can ensure
that~$\card{S}\le\card{A}\cdot\aleph_0$, see~\cite[Chapter~3]{Hodges93}.

\subsubsection*{The construction}

We well-order $\calC$ in type~$\omega_1$ and repeatedly apply the
L\"owenheim-Skolem theorem to obtain an increasing transfinite sequence
$\langle\calC_\alpha:\alpha\in\omega_1\rangle$ of elementary
sublattices of~$\calC$, such that
$\calC_\alpha=\bigcup_{\beta<\alpha}\calC_\beta$ for limit~$\alpha$.

As noted above elementarity guarantees that $\calC_0$ is a normal distributive
lattice and so its Wallman representation~$K_0$ is a metric continuum.
We may therefore apply Lemma~\ref{lemma.first.step} and obtain a map
$\psi_0:\calC_0\to\Hyper{\halfline}$ such that
\begin{itemize}
\item $\psi_0(\emptyset)=\emptyset$, \ $\psi_0(K)=\halfline$ and
      if $F\neq\emptyset$ then $\psi_0(F)$ is not compact;
\item $\psi_0(F\cup G)=\psi_0(F)\cup\psi_0(G)$;
\item if $F_1\cap\cdots\cap F_n=\emptyset$
      then $\psi_0(F_1)\cap\cdots\cap\psi_0(F_n)$ is compact; and
\item $\naturals\subseteq\psi_0\bigl(\{x\}\bigr)$, where $x$ is the point
      fixed above.
\end{itemize}

We shall construct for each $\alpha$ a map
$\psi_\alpha:\calC_\alpha\to\Hyper{\halfline}$ with the desired properties
and such that $\psi_\beta=\psi_\alpha\restr\calC_\beta$
whenever $\beta<\alpha$.
We already have $\psi_0$ and at limit stages we simply take unions
so we only have to show how to construct~$\psi_{\alpha+1}$
from~$\psi_\alpha$.

To this end enumerate $\calC_\alpha$ as $\{c_n:n\in\omega\}$ and
$\calC_{\alpha+1}\setminus\calC_\alpha$ as $\{d_n:n\in\omega\}$.
Now for each~$n$ we can find, by elementarity, an $n$-tuple
$\langle c^n_i:i<n\rangle$ of nonempty elements of\/~$\calC_\alpha$
such that
\begin{itemize}
\item for all $i,j<n$: \ $c^n_i\cup c_j=K$ \iff/ $d_i\cup c_j=K$;
\item for all $i,j<n$: \ $c^n_i\cup c^n_j=K$ \iff/ $d_i\cup d_j=K$;
\item for all $i<n$: \ $x\in c^n_i$ \iff/ $x\in d_i$; and
\item for all $f,g\subseteq n$:
      \ $\bigcap_{i\in f}c^n_i\cap\bigcap_{j\in g}c_j=\emptyset$ \iff/
        $\bigcap_{i\in f}d_i\cap\bigcap_{j\in g}c_j=\emptyset$.
\end{itemize}
To apply elementarity simply set up the right system of equations:
add $x_i\cup c_j=K$ if $d_i\cup c_j=K$ and add $x_i\cup c_j\neq K$ otherwise;
likewise add $x_i\cup x_j=K$ or $x_i\cup x_j\neq K$ etcetera, whenever
appropriate.

Given this $n$-tuple find $m_n\in\naturals$ such that
\begin{itemize}
\item if $c^n_i\cup c_j=K$
      then $[m_n,\infty)\subseteq\psi_\alpha(c^n_i)\cup\psi_\alpha(c_j)$;
\item if $c^n_i\cup c^n_j=K$
      then $[m_n,\infty)\subseteq\psi_\alpha(c^n_i)\cup\psi_\alpha(c^n_j)$;
\item if $\bigcap_{i\in f}c^n_i\cap\bigcap_{j\in g}c_j=\emptyset$
      then $\bigcap_{i\in f}\psi_\alpha(c^n_i)\cap
            \bigcap_{j\in g}\psi_\alpha(c_j)\subseteq[0,m_n)$;
\item if $x\in c^n_i$
      then $\naturals\setminus\psi_\alpha(c^n_i)\subseteq[0,m_n)$; and
\item if $x\notin c^n_i$
      then $\naturals\cap\psi_\alpha(c^n_i)\subseteq[0,m_n)$.
\end{itemize}
We can and will assume in addition that $m_{n+1}>m_n$ for all~$n$ and also
that $\psi_\alpha(c^n_i)\cap[m_n,m_{n+1})\neq\emptyset$ for all~$i$;
this is possible because $\psi_\alpha(c^n_i)$ is not compact.

Now we define, for all $i$,
$$
\psi_{\alpha+1}(d_i)=\bigcup_{n>i}\psi_\alpha(c^n_i)\cap[m_n,m_{n+1}).
$$
We show that $\psi_{\alpha+1}$ has all the required properties.

The first thing to check is that $\psi_{\alpha+1}(d_i)$ is indeed closed.
Note that only the points $m_{n+1}$ with~$n>i$ can possibly be in
$\cl\psi_{\alpha+1}(d_i)$ but not in $\psi_{\alpha+1}(d_i)$.
But if $m_{n+1}\in\cl\psi_{\alpha+1}(d_i)$
then $m_{n+1}\in\cl\bigl(\psi_{\alpha}(c^n_i)\cap[m_n,m_{n+1})\bigr)$
or $m_{n+1}\in\psi_\alpha(c^{n+1}_i)$;
in the latter case we are done and in the former case we have
$m_{n+1}\in\psi_{\alpha}(c^n_i)$ and hence $x\in c^n_i$, by the choice
of\/~$m_n$, and so $x\in d_i$.
It then follows that $x\in c^{n+1}_i$ and
so again $m_{n+1}\in\psi_\alpha(c^{n+1}_i)\subseteq\psi_{\alpha+1}(d_i)$.

It is now a routine matter to check that $\psi_{\alpha+1}$ is as required.
For example, if $d_i\cup d_j=K$ then
$[m_n,\infty)\subseteq\psi_{\alpha+1}(d_i)\cup\psi_{\alpha+1}(d_j)$,
where $n=\max(i,j)+1$ and
if $\bigcap_{i\in f}d_i\cap\bigcap_{j\in g}c_j=\emptyset$
then $\bigcap_{i\in f}\psi_{\alpha+1}(d_i)\cap
      \bigcap_{j\in g}\psi_\alpha(c_j)\subseteq[0,m_n)$,
where $n=\max(f\cup g)+1$.

In the end $\psi=\bigcup_{\alpha<\omega_1}\psi_\alpha$ is the map that we
want.

\section{A topological proof}
\label{sec.topology}

In this section we show how to modify the argument of
B\l aszczyk and Szyma\'nski from~\cite{BlaszczykSzymanski80}
to obtain a topological proof of Theorem~\ref{main.thm}.
Although this may seem to be overdoing things somewhat, we feel
it is worthwhile to include such an argument because it illustrates
our point that arguments using elementarity are necessary in the
context of connected spaces.

Indeed, consider the following statement; it seems at first sight a
reasonable property that $\Hstar$ should share with~$\Nstar$.

\begin{falselemma}\label{wishful.thinking}
If $f$ maps\/~$\Hstar$ onto the metric continuum~$X$ and if $g$ maps
the metric continuum~$Y$ onto~$X$ then there is a map~$h$ from~$\Hstar$
onto~$Y$ such that~$f=g\circ h$.
\end{falselemma}

That this lemma is false is witnessed by the following example.

\begin{example}\label{example.wish}
Define $f:\halfline\to S^1$ and $g:\unitint\to S^1$ by
$f(t)=g(t)=e^{2\pi it}$.
Now the Tietze-Urysohn theorem implies that any map from $\Hstar$
onto~$\unitint$ is induced by a map $h:\halfline\to\unitint$ and if we want
to have~$f^*=g^*\circ h^*$ then
we must have $\lim_{t\to\infty}\bigabs{f(t)-g\bigl(h(t)\bigr)}=0$.
Choose $N\in\naturals$ such that $\bigabs{f(t)-g\bigl(h(t)\bigr)}<\frac14$
for $t\ge N$ and let $n>N$.
The image, under~$f$, of $[n-\frac14,n+\frac14]$ is the arc from~$-i$,
via~$1$, to~$i$.
Therefore the image under~$g\circ h$ of this interval does not contain~$-1$
and so its image under~$h$ does not contain~$\frac12$; however this last
image does contain points above and below~$\frac12$.
This contradicts the continuity of\/~$h$.
\end{example}

The interested reader can certainly formulate the appropriate
elementary substructure version of Lemma~\ref{wishful.thinking} by
examining the inductive step in our direct algebraic construction and the
relation between $\calC_\alpha$ and $\calC_{\alpha+1}$ exploited therein.
We choose a more direct approach.

Let $K$~be the continuum of weight~$\aleph_1$ and assume that $K$~is
a subspace of the Tychonoff cube~$\unitint^{\omega_1}$.
We shall construct, by induction, maps $f_\alpha:\halfline\to\unitint$
such that the \v{C}ech-Stone extension of the diagonal
map~$f=\diag_\alpha f_\alpha$ maps $\Hstar$ onto~$K$.

For this it suffices to make sure that the following conditions are met:
\begin{enumerate}
\item \label{cond.1}
      for every open set $U$ around~$K$ there is $n\in\naturals$ such that
      $f\bigl[[n,\infty)\bigr]\subseteq U$, and
\item \label{cond.2}
      for every open set $U$ that meets~$K$ the set of $t\in\halfline$ for
      which~$f(t)\in U$ is cofinal in~$\halfline$.
\end{enumerate}

Again we shall be using a chain of elementary structures to guide our
induction but this time we need substructures of a structure for the language
of set theory.
Such a structure is $H(\kappa)$, the set of sets that are hereditarily
of cardinality less than~$\kappa$.
This means that $x\in H(\kappa)$ \iff/ $\bigcard{\tc(x)}<\kappa$,
where
$$
\tc(x)=x\cup\bigcup x\cup\bigcup\bigcup x\cup\cdots
$$
is the \emph{transitive closure} of a set~$x$.

Within such an $H(\kappa)$ one can do a lot of set theory and this makes these
sets reasonable substitutes for the whole universe of set theory,
which itself cannot be handled as an individual --- see~\cite{Kunen80b} for
the basic properties of the~$H(\kappa)$.

For our construction we take $\kappa=(2^\cont)^+$; this $\kappa$~is
big enough to guarantee that $K$, $\halfline$ and $\unitint^{\omega_1}$
belong to~$H(\kappa)$.
We also fix a continuous chain $\langle M_\alpha:\alpha<\omega_1\rangle$
of countable elementary substructures of\/~$H(\kappa)$ such that
\begin{itemize}
\item $K$ belongs to~$M_0$, and
\item for every $\alpha$ the chain $\langle M_\beta:\beta<\alpha\rangle$
      belongs to~$M_\alpha$.
\end{itemize}
Denote $M_\alpha\cap\omega_1$ by~$\delta_\alpha$ and denote the projection
of $K$ onto the first $\delta_\alpha$~coordinates by~$K_\alpha$.

It would take us too far afield to develop all possible properties of
elementary substructures of~$H(\kappa)$ but many familiar objects belong
automatically to~$M_0$ because they can be seen as the unique solution to
an equation with no constants whatsoever or with constants already obtained
from such equations:
$\emptyset$~is the unique solution to $(\forall y\in x)(y\neq y)$;
$\omega$~is the smallest set that satisfies
$(\emptyset\in x)\land\bigl((\forall y\in x)(y\cup\{y\}\in x)\bigr)$;
$\omega_1$~is the first ordinal for which there is no injective map
into~$\omega$.

Likewise $\halfline$, $\unitint$ and $\unitint^{\omega_1}$ belong to~$M_0$.
The structure~$M_0$ also contains the base~$\calB$ for~$\unitint^{\omega_1}$
that is built using only open intervals in~$\unitint$ with rational end points,
again because this set is the unique solution to an equation
involving previously identified elements of~$M_0$; this equation is nothing
but the defining sentence for~$\calB$ written out in full in the language
of set theory with $\calB$ replaced by~$x$.

We shall construct the maps $f_\beta$ for $\beta<\delta_\alpha$, by induction
on~$\alpha$.
To facilitate this we assume that the point~$\0$ with all coordinates zero
belongs to~$K$; this entails no loss of generality
because $\unitint^{\omega_1}$~is homogenous.

We use the proof of Theorem~\ref{thm.jan} to find a map
$h:\halfline\to\unitint^{\delta_0}$ such that
conditions~\ref{cond.1} and~\ref{cond.2} are met with $K$ replaced
by~$K_0$.
We can and will assume that $h(n)=\0$ for all $n\in\naturals$ and
that $h\in M_1$.
We let $f_\beta$ be the $\beta$th coordinate of\/~$h$ of course.

Now let $\alpha\in\omega_1$ and assume that the maps $f_\beta$ have been
found for $\beta<\delta_\alpha$ such that
$e_\alpha=\diag_{\beta<\delta_\alpha}f_\beta$ belongs to $M_{\alpha+1}$ and
meets conditions~\ref{cond.1} and~\ref{cond.2} with $K$ replaced
by~$K_\alpha$, and such that $f_\beta(n)=0$ for all~$\beta$ and~$n$.

Denote the family of elements of\/~$\calB$ that are supported
in~$\delta_{\alpha+1}$ by~$\calC$ and observe
that $\calC\subseteq M_{\alpha+1}$ --- every element is defined from
finitely many elements of\/~$M_{\alpha+1}$: finitely many ordinals
and rational numbers.

Let $\{U_n:n\in\omega\}$ enumerate a decreasing neighbourhood base
for~$K_{\alpha+1}$, where each $U_n$ is a finite union of elements
of\/~$\calC$ that all meet~$K_{\alpha+1}$ and are arcwise connected.
Furthermore we enumerate the elements of\/~$\calC$ that meet~$K_{\alpha+1}$
as $\{V_n:n\in\omega\}$.

Let $n\in\omega$ and consider $U_n$ and the sets $V_i$ for $i\le n$.
The set of ordinals determined by these open sets is finite
and may be split into two parts: $F_{n,0}$, those below~$\delta_\alpha$,
and~$F_{n,1}$, the rest.
One can express the facts that $K_{\alpha+1}\subseteq U_n$ and that
each~$V_i$ meets~$K_{\alpha+1}$ in one formula that involves as constants
the elements of\/~$F_{n,1}$ and some elements of\/~$M_\alpha$, notably~$K$, the
elements of\/~$F_{n,0}$ and some rational numbers.
For example if
$V_0=\pi_{\eta_1}\inv\bigl[(q_1,q_2)\bigr]\cap
     \pi_{\eta_2}\inv\bigl[(q_3,q_4)\bigr]$ then $V_0\cap K\neq\emptyset$
abbreviates
$$
(\exists x\in K)(q_1<x_{\eta_1}<q_2\land q_3<x_{\eta_2}<q_4)
$$

By elementarity one can find, in the interval $(\max F_{n,0},\delta_\alpha)$,
a finite set of ordinals~$G_n$ that can replace $F_{n,1}$ without
changing the validity of the formula.
Using~$G_n$ we can build open sets~$U_n'$ and~$V_{n,i}$ for $i\le n$ that
are reflections of\/~$U_n$ and the~$V_i$.
In the above example this means that if, say $\eta_1\in M_\alpha$
but $\eta_2\notin M_\alpha$, there must be an~$\epsilon$ in~$M_\alpha$
(larger than~$\eta_1$) such that
$$
(\exists x\in K)(q_1<x_{\eta_1}<q_2\land q_3<x_\epsilon<q_4)
$$
If this were the case $n=0$ then we would take
$V_{0,0}=\pi_{\eta_1}\inv\bigl[(q_1,q_2)\bigr]\cap
         \pi_{\epsilon}\inv\bigl[(q_3,q_4)\bigr]$.

We can therefore find, given a natural number~$N$, two natural numbers
$k$~and~$l$ bigger than~$N$ and such that
$e_\alpha\bigl[[k,\infty)\bigr]\subseteq U_n'$
and every~$V_{n,i}$ meets $e_\alpha\bigl[[k,l]\bigr]$.

Use this to find a strictly increasing sequence $\seq k$ in~$\naturals$
such that for every~$n$ we have
$e_\alpha\bigl[[k_n,\infty)\bigr]\subseteq U_n'$
and
$V_{n,i}\cap e_\alpha\bigl[[k_n,k_{n+1}]\bigr]\neq\emptyset$
for~$i\le n$.

Observe that the sets $F_{n,1}$ are increasing and that their union
is the interval~$[\delta_\alpha,\delta_{\alpha+1})$.
So given $\beta$ in this interval the set of\/~$n$ with~$\beta\in F_{n,1}$
is a final segment of\/~$\omega$.
For each such~$n$ let $\beta_n$ be the element of\/~$G_n$ that corresponds
to~$\beta$.
Finally define
$$
f_\beta(t)=
\begin{cases}
   0 & if $t\in[k_n,k_{n+1}]$ and $\beta\notin F_{n,1}$, and\cr
   f_{\beta_n}(t) & if $t\in[k_n,k_{n+1}]$ and $\beta\in F_{n,1}$.
\end{cases}
$$
Elementarity considerations will tell us that these $f_\beta$ are as required.
We have $e_{\alpha+1}\bigl[[k_n,k_{n+1}]\bigr]\subseteq U_n$ \emph{because}
$e_\alpha\bigl[[k_n,k_{n+1}]\bigr]\subseteq U_n'$ and
$V_i\cap e_{\alpha+1}\bigl[[k_n,k_{n+1}]\bigr]\neq\emptyset$
\emph{because} $V_{n,i}\cap e_\alpha\bigl[[k_n,k_{n+1}]\bigr]\neq\emptyset$
whenever~$i\le n$.

Finally note that that the whole construction can be considered to have taken
place in~$M_{\alpha+2}$ so that the result is in~$M_{\alpha+2}$ as well.

In the end the full diagonal mapping $\diag_{\beta<\omega_1}f_\beta$
is as required.

\begin{remark}
Let us note that in all three proofs it was necessary to map the set
$\Nstar$ to a single point.
In the first proof this was done to get an $\aleph_1$-saturated structure and in
the other proofs to make certain inductions work.
\end{remark}

\section{Other universal continua}
\label{sec.more.univ.continua}

We can exhibit several other natural universal continua:
for every~$n\in\naturals$ the remainder~$(\reals^n)^*$ maps onto~$\Hstar$
and hence onto every continuum of weight~$\aleph_1$.
Thus we find, under~$\CH$, for every~$n$ a universal continuum of
weight~$\cont$ and of dimension~$n$.

There is one more natural type of continuum that we want to deal with.
These we find in the remainder of the
space~$\Miod=\omega\times\unitint$ considered above.
For every~$u\in\Nstar$ the set
$$
\Isubu=\bigcap_{U\in u}\cl(U\times\unitint)
$$
is a component of\/~$\Mstar$ --- see~\cite[Section~2]{Hart92}.
Each of our three proofs can be modified to show that every
continuum~$\Isubu$ also maps onto every continuum of weight~$\aleph_1$.

This is easiest for the model-theoretic approach.
Indeed, it is very easy to show that $\Isubu$~has an $\aleph_1$-saturated
base for its closed sets:
let $\calI$ be the family of finite unions of closed intervals in~$\unitint$
with rational end~points.
The ultrapower $\calI^\omega/u$ is naturally isomorphic to a base for the
closed subsets of~$\Isubu$.
It is well-known that such ultrapowers are $\aleph_1$-saturated
--- see \cite[Theorem~9.5.4]{Hodges93}.
The isomorphism between $\calI^\omega/u$ and a base for the closed sets
of~$\Isubu$ is obtained by sending $F\in\calI^\omega$ to
$F_u=\Isubu\cap \bigl(\bigcup_{n\in\omega}(\{n\}\times F_n)\bigr)^*$;
standard properties of the \v{C}ech-Stone compactification imply
$F_u=G_u$ \iff/ $\{n:F_n=G_n\}\in u$ and that unions and intersections
are preserved under this map.

To show how to adapt the other two proofs we use the map~$q$ from~$\Mstar$
onto~$\Hstar$ defined by $q\bigl(\langle n,x\rangle\bigr)=n+x$.
The modifications will be such that the final map will map each of the
continua~$q[\Isubu]$ onto the target continuum as well.

The first step is to adapt the proof of Theorem~\ref{thm.jan}.
This is quite straightforward:
make sure that, for every~$n$, the finite set $\{a_{n,1},\ldots,a_{n,k_n}\}$
contains the points~$a_0$, \dots,~$a_{n-1}$; it then follows easily that
$A\subseteq\beta e[\Isubu]$ and hence $\beta e[\Isubu]=M$ for every~$u$.

\subsubsection*{The direct algebraic construction}
As a result of the adaptation of Theorem~\ref{thm.jan} we know
that we can ensure that for every $F\in\calC_0$ there is a natural number~$n$
such that $\psi_0(F)$ meets $[m,m+1]$ for every~$m\ge n$.

All that is needed to do then is to choose the numbers~$m_n$ so large that
for every~$i<n$ and every~$m\ge m_n$ the set $\psi_\alpha(c_i^n)$
meets the interval~$[m,m+1]$.

\subsubsection*{The topological proof}
In this case the fact that $h$ maps every $q[\Isubu]$ onto~$K_0$ translates
to: if $V$ is an open set that meets~$K_0$ then there is a natural number~$n$
such that $V$ meets~$h\bigl[[m,m+1]\bigr]$ for every~$m\ge n$.

The modification needed now is to demand each time that $V_{n,i}$ meet
$e_\alpha\bigl[[k,k+1]\bigr]$ for all~$k>k_n$.

\smallskip

In \cite{DowHart93} it is shown that, under $\CH$, all continua~$\Isubu$
are mutually homeomorphic.
This gives us a second one-dimensional universal continuum of
weight~$\cont$.
It is different from~$\Hstar$ because $\Hstar$ is indecomposable
and $\Isubu$~is not.

\section{Extensions and limitations}
\label{sec.ext.and.lim}

In this section we find analogues of some standard results about
continuous images of\/~$\Nstar$.

First we extend our main result by showing that under Martin's Axiom
($\MA$) every continuum of weight less than~$\cont$ is a continuous image
of\/~$\Hstar$.

Next we give some limiting results.
We show that in the Cohen model the long segment of
length~$\omega_2$ is not a continuous image of\/~$\Hstar$.
This shows that an extra assumption like~$\MA$ is necessary
to push Theorem~\ref{main.thm} from~$\aleph_1$ to~$\aleph_2$.
The Cohen model result also shows that in Theorem~\ref{main.thm} we cannot
replace~$\aleph_1$ by~$\cont$.
We improve this by showing that the Proper Forcing Axiom ($\PFA$) implies
that there is a continuum of weight~$\cont$ that is not a continuous image
of\/~$\Hstar$, so not even $\MA$~is strong enough to allow us to prove
Theorem~\ref{main.thm} for~$\cont$ instead of~$\aleph_1$.

\subsection*{Martin's Axiom}
\label{subsec.MA}

We show that Martin's Axiom ($\MA$) implies that all continua
of weight less than~$\cont$ are continuous images of\/~$\Hstar$.
We need the following result from~\cite{vanDouwenPrzymusinski80}.

\begin{theorem}[$\MA$]
If $X$ is a compact space of weight less than~$\cont$ then there is
a compactification of\/~$\omega$ with $X$ as its remainder.
\end{theorem}

Using this theorem and the ideas from the proof of
Theorem~\ref{thm.jan} it is now quite straightforward to prove the following
theorem by an extra application of\/~$\MA$.

\begin{theorem}[$\MA$]\label{thm.MA}
Let $K$ be a continuum of weight less than~$\cont$.
Then $K$ is a continuous image of\/~$\Hstar$.
\end{theorem}

\begin{proof}
It suffices to construct a compactification of\/~$\halfline$ with $K$ as its
remainder.
Let $\kappa$ be the weight of\/~$K$ and let $\gammaN$ be a compactification
of\/~$\omega$ with $K$ as its remainder.
We embed $\gammaN$ into~$\unitint^\kappa$ and we note the following two facts.
\begin{enumerate}
\item If $U$ is a neighbourhood of\/~$K$ in~$\unitint^\kappa$ then
      $\omega\setminus U$~is finite; and
\item $K$~has arbitrarily small arcwise connected neighbourhoods.
\end{enumerate}
The first fact is immediate; the second follows by noting that every basic
open set in~$\unitint^\kappa$ is arcwise connected and that every
neighbourhood of\/~$K$ contains a neighbourhood that is both a finite union
of basic open sets and connected.

There is a map $h:\halfline\to\unitint^\kappa$ that almost gives the desired
compactification of\/~$\halfline$: the piecewise linear map that connects
the natural numbers in their natural order.
This map is defined coordinatewise by
$h_\alpha(n+t)=(1-t)n_\alpha+t(n+1)_\alpha$.
Here $n_\alpha$ denotes the $\alpha$th coordinate of the point of\/~$\gammaN$
that corresponds to~$n$.
We need to rectify a few problems:
\ 1)~$h$ need not be one-to-one,
\ 2)~the intersection $h[\halfline]\cap K$ may be nonempty, and
\ 3)~$K$ may be a proper subset
     of\/~$\bigcap_{n\in\omega}\cl h\bigl[[n,\infty)\bigr]$.

The first two may be solved by introducing one more coordinate:
the map $t\mapsto\bigl<2^{-t},h(t)\bigr>$ from $\halfline$ to
$\unitint\times\unitint^\kappa$ does not suffer from problems~1 and~2.

To remedy~3 we introduce the following partial order~$\Poset$:
an element of\/~$\Poset$ is of the form
$p=\langle n_p,f_p,\epsilon_p,F_p,U_p\rangle$, where
$n_p\in\omega$, \ $U_p$~is an arcwise connected
neighbourhood of\/~$K$ that contains~$\omega\setminus n_p$,
 \ $f_p:[0,n_p]\to\unitint^\kappa$ is such that $f_p(n)=n$ for $n\le n_p$
--- so in particular~$f_p(n_p)\in U_p$,
\ $\epsilon_p>0$ and $F_p\in[\kappa]^{<\omega}$.

The ordering is defined as follows:
$p\le q$ \iff/
$n_p\ge n_q$, \ $\epsilon_p\le\epsilon_q$, \ $U_p\subseteq U_q$,
\ $F_p\supseteq F_q$, \ $f_p\bigl[[n_q,n_p]\bigr]\subseteq U_q$
and for all $x\in[0,n_q]$ and for all~$\alpha\in F_q$
we have $\bigl|\pi_\alpha(f_p(x))-\pi_\alpha(f_q(x))\bigr|<\epsilon_q$.

It should be clear that by making $F_p$ larger or $\epsilon_p$ and $U_p$
smaller one automatically obtains a stronger condition in~$\Poset$;
it follows that the following sets are dense in~$\Poset$:
\begin{itemize}
\item $D_{1,\alpha}=\{p:\alpha\in F_p\}$, where $\alpha\in\kappa$;
\item $D_{2,n}=\{p:\epsilon_p<2^{-n}\}$, where $n\in\omega$; and
\item $D_{3,U}=\{p:U_p\subseteq U\}$, where $U$ is a basic neighbourhood
      of\/~$K$.
\end{itemize}
We must also consider $D_{4,n}=\{p:n_p\ge n\}$, where $n\in\omega$.
To see that such a set is dense let $n$~be given.
Now if $p\in\Poset$ and $n_p<n$ then, because $U_p$ is arcwise connected and
$\{n_p,\ldots,n\}\subseteq U_p$, we can extend $f_p$ to
$g:[0,n]\to\unitint^\kappa$ with $g\bigl[[n_p,n]\bigr]\subseteq U_p$.
It follows that $\langle n,g,\epsilon_p,F_p,U_p\rangle$ is in~$D_{4,n}$
and below~$p$.

If $G$ is a filter that meets all of the aforementioned $\kappa$~many
dense sets then the net $\langle f_p\rangle_{p\in G}$ converges uniformly
on each compact subset of\/~$\halfline$ to a map~$f$, which is
continuous.
For this map~$f$ we do have
$K=\bigcap_{n\in\omega}\cl f\bigl[[n,\infty)\bigr]$.
If we add the extra coordinate to~$f$ then we get an embedding
$\bar f:t\mapsto\bigl<2^{-t},f(t)\bigr>$ of $\halfline$ into
$\unitint\times\unitint^\kappa$ that satisfies
$\cl \bar f[\halfline]=\bar f[\halfline]\cup K$.

It remains to show that $\Poset$ satisfies the countable chain condition.
So let $Q$ be an uncountable subset of\/~$\Poset$.
We may and do assume without loss of generality that there are fixed~$n$
and $\epsilon$ such that $n_p=n$ and $\epsilon_p\ge\epsilon$ for
all~$p\in Q$.
We also assume that $\{F_p:p\in Q\}$ is a $\Delta$-system with root~$F$.

With the uniform metric, the function space $C\bigl([0,n],\reals^F\bigr)$ is
separable, hence the set $\{\pi_F\circ f_p:p\in Q\}$ has a complete
accumulation point~$f$ in it.
We therefore assume, again without loss of generality, that for all
$p\in Q$ we have~$d(f,\pi_F\circ f_p)<\epsilon$.

We show that any two elements of\/~$Q$ are compatible.
Indeed, given $p$ and $q$ define $g:[0,n]\to\unitint^\kappa$ by
$$
g_\alpha(x)=\begin{cases}
            f_\alpha(x)&if $\alpha\in F$,\cr
            f_{p,\alpha}(x)&if $\alpha\in F_p\setminus F$,\cr
            f_{q,\alpha}(x)&if $\alpha\in F_q\setminus F$, and\cr
            h_\alpha(x)& otherwise.\cr
            \end{cases}
$$
The quintuple $\langle n,g,\epsilon,F_p\cup F_q,U_p\cap U_q\rangle$
is below both $p$ and~$q$.
\end{proof}

\subsection*{The Cohen Model}
\label{subsec.kunen}

In \cite[Chapter 12]{Kunen68} Kunen investigated for what
cardinals~$\kappa$ the $\sigma$-algebra~$\calR(\kappa)$ generated by the
family of rectangles equals the full power set of~$\kappa^2$, where a
\emph{rectangle} is a set of the form $X\times Y$ with
$X,Y\subseteq\kappa$.

The following lemma, from~\cite{Kunen68}, gives a convenient sufficient
condition for a set to belong to~$\calR(\kappa)$.

\begin{lemma}\label{lemma.kunen}
If $X\subseteq\kappa^2$ and if one can find, for all~$\alpha$,
subsets~$x_\alpha$ and $y_\alpha$ of~$\omega$ such that
$\langle\alpha,\beta\rangle\in X$ \iff/ $x_\alpha\cap y_\beta$ is
infinite,
then $X$ belongs to~$\calR(\kappa)$.
\end{lemma}

\begin{proof}
Define $X_n=\{\alpha:n\in x_\alpha\}$ and $Y_n=\{\alpha:n\in y_\alpha\}$.
Now observe that $\langle\alpha,\beta\rangle\in X$ \iff/
$(\forall m)(\exists n\ge m)(n\in x_\alpha\cap y_\beta)$.
Therefore $X=\bigcap_{m}\bigcup_{n\ge m} X_n\times Y_n$.
\end{proof}

Using this lemma Kunen showed that $\calR(\omega_1)=\pow(\omega_1^2)$
and that $\MA$ implies $\calR(\kappa)=\pow(\kappa^2)$ for every
$\kappa<\cont$.

Kunen complemented this by showing that in the model obtained by adding
(at least) $\aleph_2$ Cohen reals to a model of\/~$\CH$ the set
$L=\bigl\{\langle\alpha,\beta\rangle:\alpha<\beta<\omega_2\bigr\}$
does not belong to~$\calR(\omega_2)$.

\medskip
Using the results mentioned above it is nearly immediate that in the
Cohen model the ordinal space $\omega_2+1$ is not a continuous image
of~$\Nstar$; this is very likely the reason that several authors refer
to~\cite{Kunen68} as the source for this result.
Indeed, from a continuous onto mapping one immediately obtains
a sequence~$\langle x_\alpha:\alpha<\omega_2\rangle$ of subsets
of~$\omega$ such that $x_\beta\subset^* x_\alpha$ whenever
$\alpha<\beta$.
For each~$\alpha$ let $y_\alpha$ be the complement of~$x_\alpha$.
It follows that $\alpha<\beta$ \iff/ $x_\alpha\cap y_\beta$ is infinite
and so, by Lemma~\ref{lemma.kunen}, $L\in\calR(\omega_2)$.

We show that a similar result can be formulated and proved for $\Hstar$.

\begin{theorem}
In the model obtained by adding (at least)\/ $\aleph_2$ Cohen reals to
a model of\/~$\CH$ the long segment of length~$\omega_2$ is not a
continuous image of\/~$\Hstar$.
\end{theorem}

The \emph{long segment} $L_{\omega_2}$ of length~$\omega_2$ is the set
$\omega_2\times[0,1)\cup\{\omega_2\}$ ordered lexicographically.

We need the following adaptation of Lemma~\ref{lemma.kunen}.

\begin{lemma}\label{lemma.kunen.adapt}
If there is, with respect to the mod~finite order,
a strictly increasing $\omega_2$-sequence
$\langle f_\alpha:\alpha<\omega_2\rangle$
in\/~$\reals^\omega$ then $L\in\calR(\omega_2)$.
\end{lemma}

\begin{proof}
Much as in the proof of Lemma~\ref{lemma.kunen} define, for $n\in\omega$
and $q\in\rationals$, sets $X_{n,q}$ and $Y_{n,q}$ by
$X_{n,q}=\{\alpha:q>f_\alpha(n)\}$ and
$Y_{n,q}=\omega_2\setminus X_{n,q}$.
It is easily seen that
$L=\bigcap_{m}\bigcup_{n\ge m,q\in\rationals}X_{n,q}\times Y_{n,q}$.
\end{proof}

To show that in the Cohen model there is no map from $\Hstar$
onto~$L_{\omega_2}$ we show that the existence of such a map entails
the existence of sequence $\langle x_\alpha:\alpha<\omega_2\rangle$
of subsets of~$\omega$ such that $x_\beta\subset^* x_\alpha$ whenever
$\alpha<\beta$ or of a sequence $\langle f_\alpha:\alpha<\omega_2\rangle$
in~$\reals^\omega$ as in Lemma~\ref{lemma.kunen.adapt}.

So assume $f:\Hstar\to L_{\omega_2}$ is a continuous surjection.
Choose, for each $\alpha<\omega_2$, a standard open set~$U_\alpha$
in~$\betaH$ that contains $f\inv\bigl[[\alpha+1,\omega_2]\bigr]$ but whose
closure is disjoint from~$f\inv\bigl[[0,\alpha]\bigr]$.
A \emph{standard} open set is one that can be written as
$\bigcup_{n\in\omega}(a_n,b_n)$, where $a_n<b_n<a_{n+1})$ for all~$n$
and $\lim_{n\to\infty}a_n=\infty$.

Let $\langle I_n:n\in\omega\rangle$ be the sequence of intervals
that determines~$U_0$.
For each $\alpha<\omega_2$ we let
$J_\alpha=\{n: I_{0,n}\cap U_\alpha\neq\emptyset\}$.
Now either $\langle J_\alpha:\alpha<\omega_2\rangle$ is decreasing with
respect to~$\subset^*$, in which case we are done, or it becomes constant
modulo finite sets and then we may as well assume that it is constant
and that $J_\alpha=\omega$ for all~$\alpha$.

Now observe that if $\beta<\alpha$ we
have $\cl U_\alpha\cap\Hstar\subseteq U_\beta$
so there must be an~$m$ such that
$\cl U_\alpha\setminus U_\beta\subseteq[0,m]$.
We can therefore define functions $f_\alpha:\omega\to\reals$ by
$f_\alpha(n)=\inf U_\alpha\cap I_n$.
The sequence $\langle f_\alpha:\alpha<\omega_2\rangle$ is strictly increasing
with respect to~$<^*$ and we are done.

\subsection*{The Proper Forcing Axiom}

The obvious question whether $\MA$ is strong enough to imply that
every compact space of weight~$\cont$ is a continuous image of\/~$\Nstar$
has a negative answer; see for example the surveys by
Baumgartner~\cite{Baumgartner84} and Van~Mill~\cite{vanMill84}.

We show, by obvious modifications of the arguments from~\cite{Baumgartner84},
that $\PFA$~implies there is a continuum of weight~$\cont$ that is not a
continuous image of\/~$\Hstar$.

As $\PFA$ implies $\MA+{\cont=\aleph_2}$
(this was proved by Todor\v{c}evi\'c, see~\cite{Bekkali91})
we know that there is an
$\aleph_2$-saturated ultrafilter~$u$ on~$\omega$, see~\cite{EllentuckRucker72}.
Such an ultrafilter has the property that every ultrapower over it
of any countable structure is~$\aleph_2$-saturated.
We apply this to the linearly ordered set~$\rationals$ of rational numbers
and find an embedding of the lexicographically ordered tree~$T=2^{<\omega_2}$
into~$\rationals^\omega/u$.
One can map $\unitint^\omega_u$ and hence~$\rationals^\omega_u$ in a
one-to-one fashion into~$\Isubu$: send $x\in\unitint^\omega$ to the unique
point~$x_u$ in the intersection of
$\bigl\{\langle n,x_n\rangle:n\in\omega\bigr\}^*$ and~$\Isubu$.
It is explained in~\cite{Hart92} that this map is order-preserving with respect
to a natural linear quasi-order~$\prec$ in~$\Isubu$; it suffices
to know that $x_u\prec y_u$ \iff/ $\{n:x_n<y_n\}\in u$.
All this gives us an embedding of the tree~$T$ into~$\Isubu$.
We show that this standard continuum is not a continuous image of\/~$\Hstar$.

Indeed, assume there is a continuous map~$h$ from $\Hstar$ onto~$\Isubu$.
Let $x\in2^{\omega_2}$ be such that both~$x\inv(0)$ and~$x\inv(1)$
are cofinal in~$\omega_2$; then $x$~determines a
$\langle\omega_2,\omega_2^*\rangle$-gap
$\bigl<l(x),r(x)\bigr>$ in~$T$, where
$$
\begin{eqalign}
l(x)&=\bigl\{s\in T:(\exists\alpha)(\dom s=\alpha+1,
           s\restr\alpha=x\restr\alpha, s(\alpha)=0, x(\alpha)=1)\bigr\}\\
r(x)&=\bigl\{s\in T:(\exists\alpha)(\dom s=\alpha+1,
           s\restr\alpha=x\restr\alpha, x(\alpha)=0, s(\alpha)=1)\bigr\}.
\end{eqalign}
$$
Via the embedding of $T$ into~$\Isubu$ and the map~$h$ we can associate
with this gap two sequences $\{U_s:s\in l(x)\}$ and $\{V_s:s\in r(x)\}$
of standard open sets with the properties that
\begin{enumerate}
\item if $s\in l(x)$ and $t\in r(x)$ then $U_s\cap V_t=\emptyset$;
\item if $s<t$ in $l(x)$ then $\cl U_s\subseteq U_t$; and
\item if $t<s$ in $r(x)$ then $\cl V_s\subseteq V_t$.
\end{enumerate}
Now the proof, given in~\cite{Baumgartner84}, that in the Boolean
algebra~$\baB$ (the Boolean algebra of clopen subsets of~$\Nstar$) there are no
$\langle\omega_2,\omega_2^*\rangle$-gaps
(assuming~$\PFA$) works \emph{mutadis mutandis} to show that no such gaps
exist in the lattice~$\calL$ either.
It follows that there must be~$L_x\in\calL$ such that
$U_s\subseteq L_x\subseteq \Hstar\setminus V_t$ for all $s\in l(x)$
and~$t\in r(x)$.

The assignment $x\mapsto L_x$ is clearly one-to-one;
this is impossible however as there are $2^{\aleph_2}$ possible~$x$'s and
only $\aleph_2$ possible~$L_x$'s.

\subsection*{Yet another universal continuum}

In the previous subsection we used a $\cont$-saturated ultrafilter~$u$
to find a continuum of weight~$\cont$, to wit~$\Isubu$, that is not
a continuous image of\/~$\Hstar$.
The methods used in this paper allows us to show that~$\Isubu$ is in fact
a universal continuum of weight~$\cont$.

As in Section~\ref{sec.more.univ.continua} consider the family $\calI$
of all finite unions of closed intervals in~$\unitint$ with rational
end~points.
We have seen that the ultrapower~$\calI_u=\calI^\omega/u$ is a base for the
closed sets of\/~$\Isubu$; this ultrapower is $\cont$-saturated and hence
$\cont^+$-universal.
But now the reasoning leading up to Theorem~\ref{thm.alan} can be used
to show that $\Isubu$~maps onto every continuum of weight~$\cont$ or less.

\section{Universal continua need not exist}

The referee suggested that we address the very natural question whether there
exists, in $\ZFC$, a universal compact space of weight~$\cont$.

We have seen that from a $\cont$-saturated ultrafilter~$u$ we can construct a
universal continuum of weight~$\cont$.
In a similar fashion one construct a universal compact Hausdorff space of
weight~$\cont$: take the ultrapower of the Boolean algebra~$B$ of clopen subsets
of the Cantor~set by the ultrafilter~$u$; the resulting Boolean algebra is
$\cont$-saturated and hence $\cont^+$-universal.
By Stone's duality this means that its Stone space maps onto every compact
Hausdorff space of weight~$\cont$.
Unfortunately the existence of $\cont$-saturated ultrafilters is equivalent
to the conjunction of $\MA_\mathrm{countable}$ and $2^{<\cont}=\cont$ ---
see~\cite{FremlinNyikos89} --- so $\cont$-saturated ultrafilters do not
provide the definitive answer to the referee's question.
In fact the answer to the question is negative: we shall exhibit a model
of~$\ZFC$ in which there is no compact space of weight~$\cont$ that maps onto
every continuum of weight~$\cont$.
This shows that there need not be a universal compact space of weight~$\cont$
nor a universal continuum of weight~$\cont$.

\subsection*{The partial order}

We shall employ the notation of Kunen's book~\cite{Kunen80b}, except that we
shall use dots over symbols to indicate names.
We assume that our universe~$V$ satisfies~$\GCH$ and force with the partial
order $\Poset=\Fn(\omega_3\times\omega_1,2,\omega_1)\times\Fn(\omega_2,2)$.
We shall use the following properties of~$\Poset$:
\begin{enumerate}
\item $\Poset$~preserves cofinalities and hence cardinals;
\item in the generic extension $V[G]$ we have $2^{\aleph_1}=\aleph_3$ and
      $\cont=\aleph_2$;
\item $\Poset$ satisfies the $\aleph_2$-cc and
\item $G=G_1\times G_2$, where $G_1$~is
      $V$-generic on~$\Fn(\omega_3\times\omega_1,2,\omega_1)$ and $G_2$~is
      $V[G_1]$-generic on~$\Fn(\omega_2,2)$.
\end{enumerate}
The reader should consult \cite{Kunen80b}, in particular Sections~VII.6 and
VIII.4, for proofs.

\subsection*{Preparations}

Now assume that in $V[G]$ there is a universal compact space~$X$ of
weight~$\cont$; by Alexandroff's theorem, from~\cite{Alexandroff36}, there
is a compact zero-dimensional space~$Y$ of weight~$\cont$ that maps onto~$X$.
This means that we can assume that our universal space is zero-dimensional
and hence that it is the Stone space of a Boolean algebra of
cardinality~$\cont$.
We take $\omega_2$ as the underlying set of this Boolean algebra and denote
its partial order by~$\bale$; for convenience we assume that the ordinal~$0$
is the~$0$ of the Boolean algebra.

Back in $V$ we can take, in the terminology of~\cite{Kunen80b}, a nice
name~$\baledot$ for~$\bale$ and assume that all of~$\Poset$ forces
``$\baledot$~determines a Boolean algebra on~$\omega_2$ whose Stone space is
  a universal compact space of weight~$\aleph_2$''.
Because of the $\aleph_2$-cc we can find an $\aleph_2$-sized subset~$I$
of~$\omega_3$ such that all conditions involved in~$\baledot$ have the
domain of their first coordinates contained in $I\times\omega_1$; without loss
of generality and for notational convenience we take $I=\omega_2$.

\subsection*{A nonimage}

We first show how to find a zero-dimensional compact space of weight~$\cont$
that is not a continuous image of our space~$X$; in the next subsection we show
how to modify the construction so as to obtain a continuum of weight~$\cont$.

Consider the tree $T=2^{<\omega_1}$ from~$V$.
In $V[G]$ we consider the following branches of~$T$: for each~$\alpha<\omega_2$
define $B_\alpha$ by $B_\alpha(\xi)=(\bigcup G_1)(\omega_2+\alpha,\xi)$.
Put $\calB=\{B_\alpha:\alpha\in\omega_2\}$ and consider the tree~$T\cup\calB$.
We turn this tree upside-down to make it generate a Boolean algebra~$\baT$ of
cardinality~$\cont=\aleph_2$.
We shall show that this Boolean algebra cannot be embedded into
$\orpr{\omega_2}\bale$ or, dually, that its Stone space~$Y$ is not a continuous
image of~$X$.
So we assume $\phi:T\cup\calB\to\omega_2$ is the restriction of an embedding
of~$\baT$ into $\orpr{\omega_2}\bale$ and proceed to reach a contradiction.

Back in~$V$ we take a nice name~$\phidot$ for~$\phi$.
We apply the $\aleph_2$-cc once more to find a set $J$ of size~$\aleph_1$ such
that for all~$\orpr{p}{q}$ involved in~$\phidot\restr T$ we have
$\dom(p)\subseteq J\times\omega_1$.
Put $S=\omega_2\cup J$ and set $R=\omega_3\setminus S$; apply the product lemma
to write $G_1=G_s\times G_r$, where $G_s$~is $V$-generic
on~$\Fn(S\times\omega_1,2,\omega_1)$ and $G_r$~is $V[G_s]$-generic
on~$\Fn(R\times\omega_1,2,\omega_1)$.
In fact, the results from \cite[Section~VIII.1]{Kunen80b} allow us to conclude
that $G_r$~is $V[G_s\times G_1]$-generic on~$\Fn(R\times\omega_1,2,\omega_1)$.

We reach our final contradiction by taking
$\delta\in[\omega_2,\omega_2\cdot2)\cap R$
and proving that $B_\delta\in V[G_s\times G_1]$.
Indeed: let $\gamma=\phi(B_\delta)$ and observe that $B_\delta$~is definable
from $\gamma$, $\phi\restr T$ and $\bale$:
$$
B_\delta=\bigl\{s\in T:\gamma\bale \phi(s)\bigr\},
$$
this is because $T$~is a tree and $0\prec\gamma$.
The three parameters belong to~$V[G_s\times G_1]$, hence so does~$B_\delta$.

\subsection*{A continuum}

We now show how to modify the nonimage from the previous subsection so as to get
a continuum; basically we replace~$2$ with~$\Z$, the set of integers.
So we force with $\Fn(\omega_3\times\omega_2,\Z,\omega_1)\times\Fn(\omega_2,2)$
and consider the tree $T=\Z^{<\omega_1}$.
We give $T$ the lexicographic order:
\begin{itemize}
\item if $s\restr\alpha=t\restr\alpha$ and $s(\alpha)<t(\alpha)$
      then $s\lelex t$; and
\item if $s=t\restr\alpha$ then $t\lelex s$ if $t(\alpha)<0$
      and $s\lelex t$ if $t(\alpha)\ge 0$.
\end{itemize}
Observe that $\lelex$ is a dense linear order on~$T$.

Now, as before, assume $\bale$ is a Boolean partial order on~$\omega_2$ and
let $\baledot$ be a nice name for it, which we assume to be determined by
$\Fn(\omega_2\times\omega_1,\Z,\omega_1)\times\Fn(\omega_2,2)$.
We use $G_1$ to define branches $\{B_\alpha:\alpha<\omega_2\}$ of~$T$ and
use each of these branches to insert a copy of~$\rationals$
in~$\orpr{T}{\lelex}$.
Our continuum~$K$ is the Dedekind completion of this expanded linearly
ordered set.

We show that $K$ is not a continuous image of the Stone space, call
it~$X$, of our Boolean algebra.
So assume there is a continuous surjection~$h$ of~$X$ onto~$K$.
Every $s\in T\cup\calB$ determines a closed interval~$I_s$ in~$K$, namely
the closure of $\{t:s\subseteq t\}$.
Because we used $\Z$ rather than~$2$ we know that $I_t$ is contained in the
interior of~$I_s$ whenever $s\subset t$.
Therefore, if $s$~is on a successor level, say with immediate predecessor~$s^-$,
we can take~$\phi(s)\in\omega_2$ so that
$h\inv[I_s]\subseteq \phi(s)\subseteq h\inv[I_{s^-}]$ (we identify $\phi(s)$
with the clopen subset of~$X$ that it represents).
We also choose, for each~$\alpha<\omega_2$, one element~$\gamma_\alpha$ whose
nonempty clopen set is contained in~$h\inv[I_{B_\alpha}]$.

The contradiction is reached exactly as before: one can recover $B_\delta$
from $\gamma_\delta$, $\phi$, $T$ and $\bale$ using almost the same formula:
$$
B_\delta=\Bigl\{s\in T:(\exists t\in T)\bigl(s\subseteq t\text{ and }
             \gamma_\delta\bale\phi(t)\bigr)\Bigr\}.
$$

\section{Some questions}

Of course every result about continuous images of\/~$\Nstar$ suggests
the possibility of a corresponding one about continuous images of\/~$\Hstar$.
The more obvious questions~are.

\begin{question}[compare \cite{Przymusinski82}]
Is every perfectly normal continuum a continuous image of\/~$\Hstar$?
\end{question}

\begin{question}
Is it consistent with $\lnot\CH$ (or better still with $\MA+\lnot\CH$)
that every continuum of weight~$\cont$ is a continuous image of\/~$\Hstar$?
\end{question}

A more interesting question, raised by G.~D. Faulkner and also suggested by
the proofs of Theorems~\ref{thm.jan} and~\ref{thm.MA}, is.

\begin{question}
Is every continuum that is a continuous image of\/~$\Nstar$ also a continuous
image of\/~$\Hstar$?
\end{question}

A still more interesting problem is to generalize Parovi\v{c}enko's
characterization of\/~$\Nstar$ to the connected case.

\begin{problem}
Find, assuming $\CH$, topological characterizations of the continua
$\Hstar$ and $\Isubu$.
\end{problem}

This problem could have quite interesting consequences since such a
characterization would use some specific base for the closed sets
of~$\Hstar$ and, as we see by the methods in this paper,
it is likely that it is first necessary to find a natural description of the
first-order theory of the lattice of discrete unions of intervals,
or of some other natural base for the closed sets of~$\halfline$.


\begin{thebibliography}{10}

\bibitem{AartsvanEmdeBoas67}
J.~M. Aarts and P.~van Emde~Boas, \emph{Continua as remainders in compact
  extensions}, Nieuw Archief voor Wiskunde (3) \textbf{15} (1967), 34--37.

\bibitem{Alexandroff27}
P.~S. Alexandroff, \emph{{\"Uber} stetige {Abbildungen} kompakter {R\"aume}},
  Mathematische Annalen \textbf{96} (1927), 555--571.

\bibitem{Alexandroff36}
\bysame, \emph{Zur {Theorie} der topologischen {R{\"a}ume}}, Comptes Rendus
  (Doklady) de l'Acad{\'e}mie des Sciences de l'URSS \textbf{11} (1936),
  55--58.

\bibitem{Baumgartner84}
James~E. Baumgartner, \emph{Applications of the {Proper} {Forcing} {Axiom}}, In
  Kunen and Vaughan \cite{KunenVaughan84}, pp.~913--960.

\bibitem{Bekkali91}
Mohamed Bekkali, \emph{Topics in set theory}, Lecture Notes in Mathematics, no.
  1476, Springer-Verlag, Berlin~etc., 1991.

\bibitem{BlaszczykSzymanski80}
A.~B{\l}aszczyk and A.~Szyma{\'n}ski, \emph{Concerning {P}arovi{\v{c}}enko's
  theorem}, Bulletin de L'Academie Polonaise des Sciences S\'{e}rie des
  sciences math\'ematiques \textbf{28} (1980), 311--314.

\bibitem{vanDouwen81d}
Eric~K. van Douwen, \emph{Special bases for compact metrizable spaces},
  Fundamenta Mathematicae \textbf{111} (1981d), 201--209.

\bibitem{vanDouwenPrzymusinski80}
Eric~K. van Douwen and T.~C. Przymusi\'nski, \emph{Separable extensions of
  first-countable spaces}, Fundamenta Mathematicae \textbf{105} (1980),
  147--158.

\bibitem{DowHart93}
Alan Dow and Klaas~Pieter Hart, \emph{{{\v{C}}ech-Stone} remainders of spaces
  that look like~$[0,\infty)$}, Acta Universitatis Carolinae---Mathematica et
  Physica \textbf{34} (1993), no.~2, 31--39, published in~1994.

\bibitem{EllentuckRucker72}
Erik Ellentuck and R.~V.~B. Rucker, \emph{Martin's {Axiom} and saturated
  models}, Proceedings of the American Mathematical Society \textbf{34} (1972),
  243--249.

\bibitem{FremlinNyikos89}
D.~H. Fremlin and P.~J. Nyikos, \emph{Saturating ultrafilters on~$\naturals$},
  Journal of Symbolic Logic \textbf{54} (1989), 708--718.

\bibitem{Hart92}
Klaas~Pieter Hart, \emph{The {{\v C}ech-Stone} compactification of the {Real
  Line}}, Recent Progress in General Topology (Miroslav Hu{\v{s}}ek and Jan van
  Mill, eds.), North-Holland, Amsterdam, 1992, pp.~317--352.

\bibitem{Hausdorff27}
Felix Hausdorff, \emph{Mengenlehre. {3.~Auflage}}, G{\"o}schens
  Lehrb{\"u}cherei, no.~7, De Gruyter, Berlin and Leipzig, 1935, English
  Translation: \textsl{Set Theory}, Chelsea Publications Co. New York, 1957.

\bibitem{Hodges93}
Wilfrid Hodges, \emph{Model theory}, Encyclopedia of mathematics and its
  applications, no.~42, Cambridge University Press, Cambridge, 1993.

\bibitem{JonssonOlin68}
Bjarni J\'{o}nsson and Philip Olin, \emph{Almost direct products and
  saturation}, Compositio Mathematica \textbf{20} (1968), 125--132.

\bibitem{Kunen68}
K.~Kunen, \emph{Inaccessibility properties of cardinals}, Ph.D. thesis,
  Stanford University, 1968.

\bibitem{Kunen80b}
\bysame, \emph{Set theory. an introduction to independence proofs}, Studies in
  Logic and the foundations of mathematics, no. 102, North-Holland, Amsterdam,
  1980b.

\bibitem{KunenVaughan84}
Kenneth Kunen and Jerry~E. Vaughan (eds.), \emph{Handbook of set theoretic
  topology}, North-Holland, Amsterdam, 1984.

\bibitem{Kuratowski66}
K.~Kuratowski, \emph{Topology~{I}}, PWN---Polish Scientific Publishers and
  Academic Press, Warszawa and New York, 1966.

\bibitem{vanMill84}
Jan van Mill, \emph{An {Introduction} to $\beta\omega$}, In Kunen and Vaughan
  \cite{KunenVaughan84}, pp.~503--568.

\bibitem{Parovicenko63}
I.~I. Parovi{\v{c}}enko, \emph{A universal bicompact of weight $\aleph$},
  Soviet Mathematics Doklady \textbf{4} (1963), 592--595, Russian original:
  {\cyr \emph{Ob odnom universal{\cprime}nom bikompakte vesa~$\aleph$}, Doklady
  Akademii Nauk SSSR \textbf{150} (1963) 36--39}.

\bibitem{Przymusinski82}
T.~C. Przymusi\'nski, \emph{Perfectly normal compact spaces are continuous
  images of $\beta\naturals-\naturals$}, Proceedings of the American
  Mathematical Society \textbf{86} (1982), 541--544.

\bibitem{StoneMH37a}
Marshall~H. Stone, \emph{Applications of the theory of {Boolean} rings to
  general topology}, Transactions of the American Mathematical Society
  \textbf{41} (1937{a}), 375--481.


\bibitem{Wallman38}
H.~Wallman, \emph{Lattices and topological spaces}, Annals of Mathematics
  \textbf{39} (1938), 112--126.

\bibitem{Waraszkiewicz34}
Z.~Waraszkiewicz, \emph{Sur un probl\`eme de {M. H. Hahn}}, Fundamenta
  Mathematicae \textbf{22} (1934), 180--205.

\end{thebibliography}

\providecommand{\bysame}{\leavevmode\hbox to3em{\hrulefill}\thinspace}

\end{document}